\documentclass{amsart}

\usepackage{amsmath,amssymb,amsthm}
\usepackage{graphicx,tikz}

\newtheorem*{thm}{Theorem}

\newtheorem*{corollary}{Corollary}

\theoremstyle{definition}

\theoremstyle{remark}

\begin{document}

\title[]{A Weighted Randomized Kaczmarz Method \\for solving linear systems}
\thanks{S.S. is supported by the NSF (DMS-1763179) and the Alfred P. Sloan Foundation.}

\author[]{Stefan Steinerberger}
\address{Department of Mathematics, University of Washington, Seattle}
\email{steinerb@uw.edu}

\begin{abstract} The Kaczmarz method for solving a linear system $Ax = b$ interprets such a system as a collection of equations $\left\langle a_i, x\right\rangle = b_i$, where $a_i$ is the $i-$th row of $A$. It then picks such an equation and corrects $x_{k+1} = x_k + \lambda a_i$ where $\lambda$ is chosen so that the $i-$th equation is satisfied. Convergence rates are difficult to establish. Strohmer \& Vershynin established that if the order of equations is chosen randomly with likelihood proportional to the size of $\|a_i\|^2_{\ell^2}$), then $\mathbb{E}~ \|x_k - x\|_{\ell^2}$ converges exponentially. We prove that if the $i-$th row is selected with likelihood proportional to $\left|\left\langle a_i, x_k \right\rangle - b_i\right|^{p}$, where $0<p<\infty$, then $\mathbb{E}~\|x_k - x\|_{\ell^2}$ converges faster than the purely random method.
As $p \rightarrow \infty$, the method de-randomizes and explains, among other things, why the maximal correction method works well. \end{abstract}

\maketitle

\section{Introduction}
Let $A \in \mathbb{R}^{m \times n}$, $m \geq n$, be a matrix with rank $n$. Let $x \in \mathbb{R}^n$ and  
$$ Ax =b,$$
how would one go about finding $x$ from $b$ and $A$?
The Kaczmarz method \cite{kac} (also known as the \textit{Algebraic Reconstruction Technique} (ART) in computer tomography \cite{gordon, her, her2, natterer} or the \textit{Projection onto Convex Sets Method} \cite{cenker, deutsch, deutsch2, gal, sezan}) is an established and well-studied iterative method for finding the solution. Denoting the rows of $A$ by $a_1, \dots, a_m \in \mathbb{R}^n$, the linear system can be written in terms of a collection of
inner products,
$$ \forall~1 \leq i \leq m: \left\langle a_i, x\right\rangle = b_i.$$
The Kaczmarz method proceeds as follows: given an approximation $x_k$, let us pick an equation, say the $i-$th equation, and construct $x_{k+1}$ from $x_k$ via
$$ x_{k+1} = x_k + \frac{b_i - \left\langle a_i, x_k\right\rangle}{\|a_i\|^2}a_i.$$
This ensures that the $i-$th equation is satisfied (though the other equations usually will not be). By cycling through the indices in a periodic fashion, we converge towards the solution $x$.
The algorithm has a geometric interpretation (see Fig. 1): $x_{k+1}$ is obtained by projecting $x_k$ onto the hyperplane given by $\left\langle a_i, x \right\rangle = b_i$. Convergence rates are difficult to obtain. Strohmer \& Vershynin \cite{strohmer} showed that by choosing the equations in a random order (with the $i-$th equation being chosen with likelihood proportional to $\|a_i\|^2$), the method converges exponentially and
$$ \mathbb{E} \left\| x_k - x \right\|_2^2 \leq \left(1 - \frac{1}{\|A\|_F^2 \| A^{-1}\|_2^2}\right)^k \|x_0 - x\|_2^2,$$
where $\|A^{-1}\|_2$ is the operator norm of the inverse and $\|A\|_F$ is the Frobenius norm.
This remarkable result has inspired a lot of subsequent work: we refer to 
\cite{bai, eldar, elf, gordon2, gower, lee, leventhal, liu, ma, moor, need, need2, need25, need3, need4, popa, tan, zhang, zouz} and references therein. It was recently shown in \cite{stein} that if $v_n \in \mathbb{R}^n$ denotes the singular vector corresponding to the smallest singular value of $A$, then
$$ \mathbb{E} \left\langle x_k - x, v_n\right\rangle =  \left(1 - \frac{1}{\|A\|_F^2 \| A^{-1}\|_2^2}\right)^k \left\langle x_0 - x, v_n\right\rangle.$$
This shows that the exponential rate obtained by Strohmer \& Vershynin \cite{strohmer} accurately captures the underlying process in the generic setting.

\begin{center}
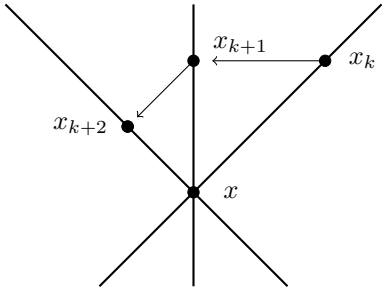
\begin{figure}[h!]
\begin{tikzpicture}[scale=2.5]
\draw[thick] (-0.5, -0.5) -- (1, 1);
\draw[thick] (0, -0.5) -- (0, 1);
\draw[thick] (0.5, -0.5) -- (-1, 1);
\filldraw (0,0) circle (0.03cm);
\node at (0.2, 0) {$x$};
\filldraw (0.7, 0.7) circle (0.03cm);
\node at (0.9, 0.7) {$x_k$};
\filldraw (0, 0.7) circle (0.03cm);
\node at (0.25, 0.78) {$x_{k+1}$};
\draw [ <-] (0.1, 0.7) -- (0.7, 0.7);
\filldraw (-0.7/2, 0.7/2) circle (0.03cm);
\filldraw (-0.35, 0.35) circle (0.03cm);
\node at (-0.6, 0.35) {$x_{k+2}$};
\draw [->] (0, 0.7) -- (-0.3, 0.4);
\end{tikzpicture}
\caption{A geometric interpretation: iterative projections onto the hyperplanes given by $\left\langle a_i, x\right\rangle = b_i$.}
\end{figure}
\end{center}

\section{Results}
\subsection{The Algorithm.} We will assume without loss of generality that the matrix is normalized to $\|a_i\|=1$ for each row $a_i$. This simplifies the description of the Randomized Kaczmarz method: by picking each equation with equal likelihood, we obtain an iterative scheme having convergence rate
$$ \mathbb{E} \left\| x_k - x \right\|_2^2 \leq \left(1 - \frac{1}{\|A\|_F^2 \cdot \|A^{-1}\|^2}\right)^k \|x_0 - x\|_2^2.$$
We were motivated by the following heuristic. Let us suppose $x_k$ is given and $\left\langle x_k, a_i \right\rangle \sim b_i$ with a very small error. Then, surely,
we should not pick the $i-$th equation in the next step: the error is already small and we gain very little. In fact, we should indeed try to correct for an equation where $\left\langle x_k, a_i\right\rangle$ is quite different from $b_i$. So a natural question is: how about we use a quantity related to the size of $\left| \left\langle x_k, a_i \right\rangle - b_i\right|$ as the likelihood of picking the $i-$th equation? The special case of picking the largest entry is quite classical: sometimes known as \textit{Motzkin's method} \cite{agmon} or the \textit{maximal correction method} by Cenker, Feichtinger, Mayer, Steier \& Strohmer \cite{cenker}. We propose a method in a similar spirit and comment on the connection to the maximal correction method in \S 2.3. \\

\textbf{A Weighted Randomized Kaczmarz Algorithm.}
\begin{quote}
Input: a matrix $A \in \mathbb{R}^{m \times n}$ with $m \geq n$ and full rank and each row normalized to $\|a_i\|=1$, a vector $b \in \mbox{Im}(A) \subseteq \mathbb{R}^m$, a parameter $0 < p < \infty$ and an initial $x_0 \in \mathbb{R}^n$. \\
Output: a sequence $(x_k)_{k=1}^{\infty}$ converging to $x=A^{-1}b$.\\
\textbf{\hspace{10pt}Pre-Computation.}\\
\vspace{-12pt}
\begin{enumerate} 
\item Compute $r_0 = Ax_0$. 
\item Compute the matrix $Q \in \mathbb{R}^{m \times m}$ given by $Q_{ij} = \left\langle a_i, a_j\right\rangle$.\\
\textbf{Generation of the sequence.}
\item Compute the vector $Ax_k - b$.
\item Use it to pick a random equation $1 \leq i \leq m$ such that 
$$ \mathbb{P}(\mbox{we choose equation}~i) =  \frac{\left| \left\langle a_i, x_k \right\rangle - b \right|^p}{\|Ax_k - b \|^p_{\ell^p}}.$$
\item Compute 
$$ \lambda = b_i - \left\langle a_i, x_k \right\rangle.$$
\item Generate the next element in the sequence
$$ x_{k+1} = x_k + \lambda a_i.$$
\item Update the right-hand side by computing 
$$ r_{k+1} = r_k + \lambda Q_i.$$
\item Update the index, $k \rightarrow k+1$, and go to (3).
\end{enumerate}
\end{quote}

We note that, a priori, the method is quite similar to Randomized Kaczmarz (which it emulates as $p\rightarrow 0$). The main difference is that the likelihood of
choosing a fixed equation is not fixed but is also being updated. We note that if it is possible to store the $n \times n$ matrix
$Q_{ij} = \left\langle a_i, a_j\right\rangle$, then this updating step is essentially as expensive as updating the sequence element 
and thus quite cheap.  If $m \gg n$, then storing $Q$ (at cost $m^2$) might be more expensive than storing $A$ (at cost $m \cdot n$),
however, the system is then also massively overdetermined and one might arguably be able to find a solution by using only a subset
of the information.

\subsection{The Result.} We can now state the main result: this algorithm is at least as efficient as the fully random version, independently of the value of $p$.
\begin{thm} Let $0 < p < \infty$, let $A$ be normalized to having the norm of each row be $\|a_i\|=1$. Then
$$ \mathbb{E} \left\| x_k - x \right\|_2^2 \leq \left(1 - \inf_{z \neq 0}  \frac{\|A z \|^{p+2}_{\ell^{p+2}}}{\|Az \|^p_{\ell^p}\|z\|^2_{2}}\right)^k \|x_0 - x\|_2^2.$$
This is at least the rate of Randomized Kaczmarz since
$$ \inf_{z \neq 0} \frac{\|A z \|^{p+2}_{\ell^{p+2}}}{\|A z \|^p_{\ell^p}\|z\|^2_{2}} \geq  \frac{1}{\|A\|_F^2  \cdot \|A^{-1}\|^2}$$
with equality if and only if the singular vector $v_n$ corresponding to the smallest singular value of $A$ has the property that $A v_n$ is a constant vector.
\end{thm}

While the rate obtained by Strohmer \& Vershynin \cite{strohmer} is sharp for
the fully random case, it is less clear whether Theorem 1 gives the optimal rate for this larger family of algorithms. We also note the monotonicity of 
$$ \frac{\|A z \|^{p+2}_{\ell^{p+2}}}{\|Az \|^p_{\ell^p}\|z\|^2_{2}} \qquad \mbox{in}~p.$$
This can be seen as follows: the quantity is invariant under multiplication of $x$ by $\lambda > 0$, so we may assume w.l.o.g. that each non-zero coordinate satisfies $|x_i| \geq 1$. Ignoring the constant factor $\|x\|_{2}^2$, we can write, denoting $\log{|x_i|} = c_i \geq 0$.
$$ \frac{\|A x \|^{p+2}_{\ell^{p+2}}}{\|Ax \|^p_{\ell^p}} = \frac{\sum_{i=1}^{n}{e^{c_i (p+2)}}}{\sum_{i=1}^{n}{e^{c_i p}}} =\sum_{i=1}^{n}e^{2 c_i}   \frac{e^{c_i p}}{\sum_{i=1}^{n}{e^{c_i p}}}.$$
 As $p$ increases, the weights get disproportionately larger for large values of $c_i$ which are in front of larger constants. Numerical experiments will also confirm that larger values of $p$ correspond to better performance. We believe that this Theorem naturally suggests a large number of questions (some of which are discussed in \S 2.3 and \S 3.4). In particular, it would be interesting to understand what happens when the linear system does not have a solution (this was done by Needell \cite{need} for the Random Kaczmarz method) or whether the assumption that $A$ has full rank can be relaxed (see Zouzias \& Freris \cite{zouz}).

 \subsection{The case $p \rightarrow \infty$.} We note that our main result is free of constants implicitly depending on $p$. It is thus not a problem to choose $p$ to be a very large number. Since we select the equations according to the size of $\left| \left\langle a_i, x_k \right\rangle - b \right|^p$, something interesting happens when $p$ gets very large: in most cases and with overwhelming likelihood, we end up simply picking an equation for which
 $$ \left| \left\langle a_i, x_k \right\rangle - b \right| = \|Ax -b \|_{\ell^{\infty}}.$$
 This is known as the \textit{maximal correction method}, we refer to  Cenker, Feichtinger, Mayer, Steier \& Strohmer \cite{cenker}. Nutini, Sepehey, Laradji, Schmidt, Koepke, Virani \cite{nutini} proved that greedy selection rules lead to superior results. We quickly discuss the implications of a statement implied by our proof. 
 \begin{corollary} Let $0 < p < \infty$. We have
 $$ \mathbb{E} \left\| x_{k+1} - x \right\|_2^2 \leq \left(1 -  \frac{\|A (x-x_k) \|^{p+2}_{\ell^{p+2}}}{\|A(x- x_k) \|^p_{\ell^p}\|x-x_k\|^2_{2}}\right) \|x_k - x\|_2^2.$$
 \end{corollary}
This implies the following dichotomy: at each iteration step
\begin{enumerate}
\item \textbf{either} $\|A(x-x_k)\|_{\ell^{p+2}}$ is quite small (and therefore $Ax_k$ already provides a good uniform approximation of $b$) 
\item \textbf{or}, if $\|A(x-x_k)\|_{\ell^{p+2}}$ is quite large, then it forces $\|x_{k+1}-x\|_{2} \ll \|x_k - x\|_{}$.
\end{enumerate}
This is perhaps one of the reasons why one would expect nice convergence properties for $p$ large: the error estimates are nicely bootstrapping themselves. Either the uniform error is small or the $\ell^2-$norm decays rapidly (eventually also forcing a small $\ell^{\infty}$ error). This can be observed numerically as well, larger values of $p$ result in approximations for which $\|A x_k -b\|_{\ell^{\infty}}$ decays more rapidly (see \S 3.3). It would be desirable to have more refined estimates for the behavior of $\|A x_k -b\|_{\ell^{\infty}}$ under these methods: do they yield good uniform approximation?

\subsection{Related Results.} This result is related to recent work of Gower, Molitor, Moorman and Needell \cite{gower2} who discuss non-uniform selection probabilities in the more general framework of Sketch-and-Project methods. In particular, \cite[\S 7.4]{gower2} discusses the convergence rate of the algorithm considered here for $p=1$. The case $p=\infty$ has been studied by a large number of people; we emphasize the results of Nutini, Sepehey, Laradji, Schmidt, Koepke, Virani \cite{nutini}.
 Haddock \& Ma \cite{haddock} recently gave a new analysis of the $p=\infty$ case, we also refer to Haddock \& Needell \cite{haddock2}. Other variations on greedy Kaczmarz have also been recently proposed by Bai \& Wu \cite{bai0, bai, bai2}, Du \& Gao \cite{du} and Li \& Zhang \cite{li}. A related approach for coordinate descent has been investigated by Li, Lu \& Wang \cite{lilu}. Most of these results are very much inspired by the original argument of Strohmer \& Vershynin \cite{strohmer} and phrase bounds in terms of classical notions of singular values and the Frobenius norm. The main insight of our paper is that choosing the probabilities in this particular algebraic way allows us to formulate better bounds by deviating from the $\ell^2-$geometry and phrasing things in $\ell^p$ (however, for $p=\infty$, something similar has already been done by Nutini, Sepehey, Laradji, Schmidt, Koepke, Virani \cite{nutini}). Note added in print: we were informed by one of the referees that our Algorithm has also appeared in a more recent preprint of Jiang, Wu \& Jiang \cite{jiang}. They argue, based on numerical evidence, that larger values of $p$ lead to better results; they then propose $p=\infty$ as a natural choice for which they provide a convergence analysis that is quite different from ours. 

\section{Numerical Examples}
\subsection{A nice matrix.} We first consider a simple example. Construct $A \in \mathbb{R}^{1000 \times 1000}$ by taking each entry to be $a_{ij} \sim \mathcal{N}(0,1)$, add $100 \cdot \mbox{Id}_{1000 \times 1000}$ and then normalize the rows. Given this matrix, we solve the problem $Ax = 0$ starting with $x_0 = (1,1,\dots, 1)$. This problem is somewhat well posed, the smallest singular value is at size $\sim 1$. The result is seen in Fig. 2, the weighted versions converge much more rapidly. Moreover, the higher the value of $p$, the more rapid the convergence.

\begin{center}
\begin{figure}[h!]
\begin{tikzpicture}[scale=1]
\node at (0,0) {\includegraphics[width=0.55\textwidth]{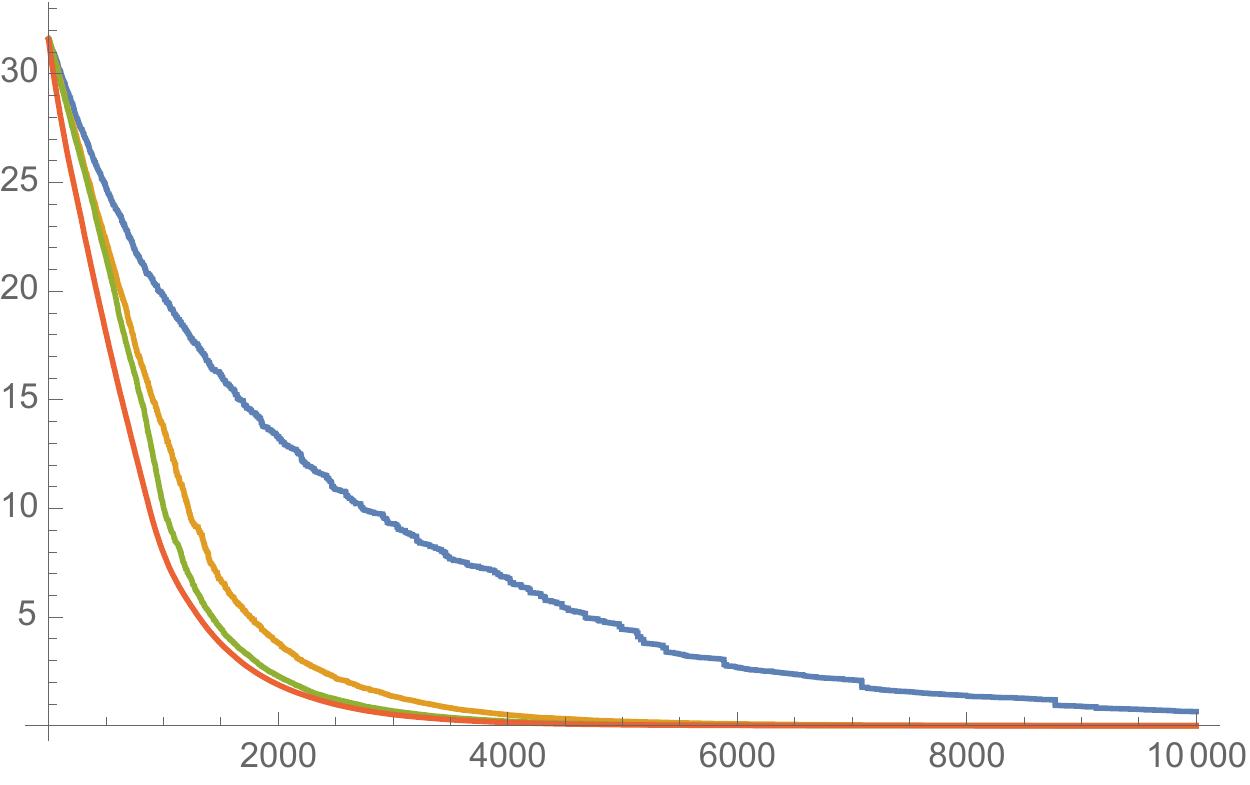}};
\node at (1,1) {$\|x_k - x\|_{\ell^2}$};
\end{tikzpicture}
\caption{ $\|x_k - x\|_{2}$ for the Randomized Kaczmarz method (blue), for $p=1$ (orange), $p=2$ (green) and $p=20$ (red).}
\end{figure}
\end{center}

\subsection{A more challenging matrix.}
We take $A \in \mathbb{R}^{1000 \times 1000}$ by taking each entry to be independently and identically distributed $a_{ij} \sim \mathcal{N}(0,1)$, then normalize each row to $\|a_i\|=1$ and solve the system $Ax = 0$ starting with $x_0 = (1,1,1\dots, 1)$. This is a hard problem, the smallest singular values being at scale $10^{-3}$.
We observe again, see Fig. 3, that larger values of $p$ correspond to a better performance. The difference in asymptotic behavior seems somewhat negligible but the initial regularization happens much more rapidly in the weighted scheme. The runtime of the various algorithms is comparable up to a factor of 2 in our implementation of Mathematica (including the computation of $Q$).

\begin{center}
\begin{figure}[h!]
\begin{tikzpicture}
\node at (0,0) {\includegraphics[width=0.55\textwidth]{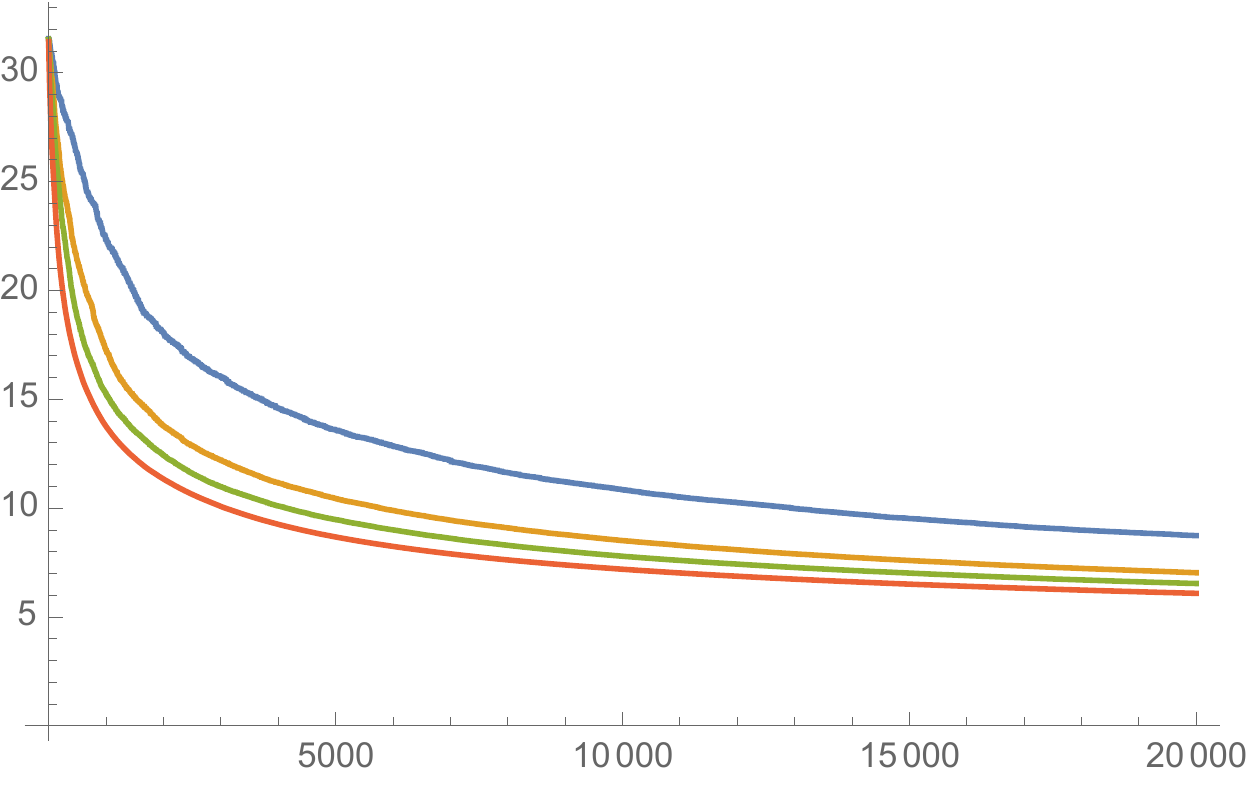}};
\node at (1,1) {$\|x_k - x\|_{\ell^2}$};
\end{tikzpicture}
\caption{$\|x_k - x\|_{2}$ for the Randomized Kaczmarz method (blue), for $p=1$ (orange), $p=2$ (green) and $p=20$ (red).}
\end{figure}
\end{center}

\subsection{Uniform Approximation.} We revisit the nice matrix from \S 3.1. above and plot the decay of the uniform error $\|Ax_k -b\|_{\ell^{\infty}}$ (see Fig. 4). That error also decays much more rapidly with larger values of $p$. We see that the Randomized Kaczmarz method, while having fairly continuous levels of decay in $\|x_k-x\|_{2}$ (see \S 3.1 and \S 3.2) does not particularly improve in terms of uniform approximation until `bad equations' are randomly chosen which then leads to a massive drop.
\begin{center}
\begin{figure}[h!]
\begin{tikzpicture}
\node at (0,0) {\includegraphics[width=0.6\textwidth]{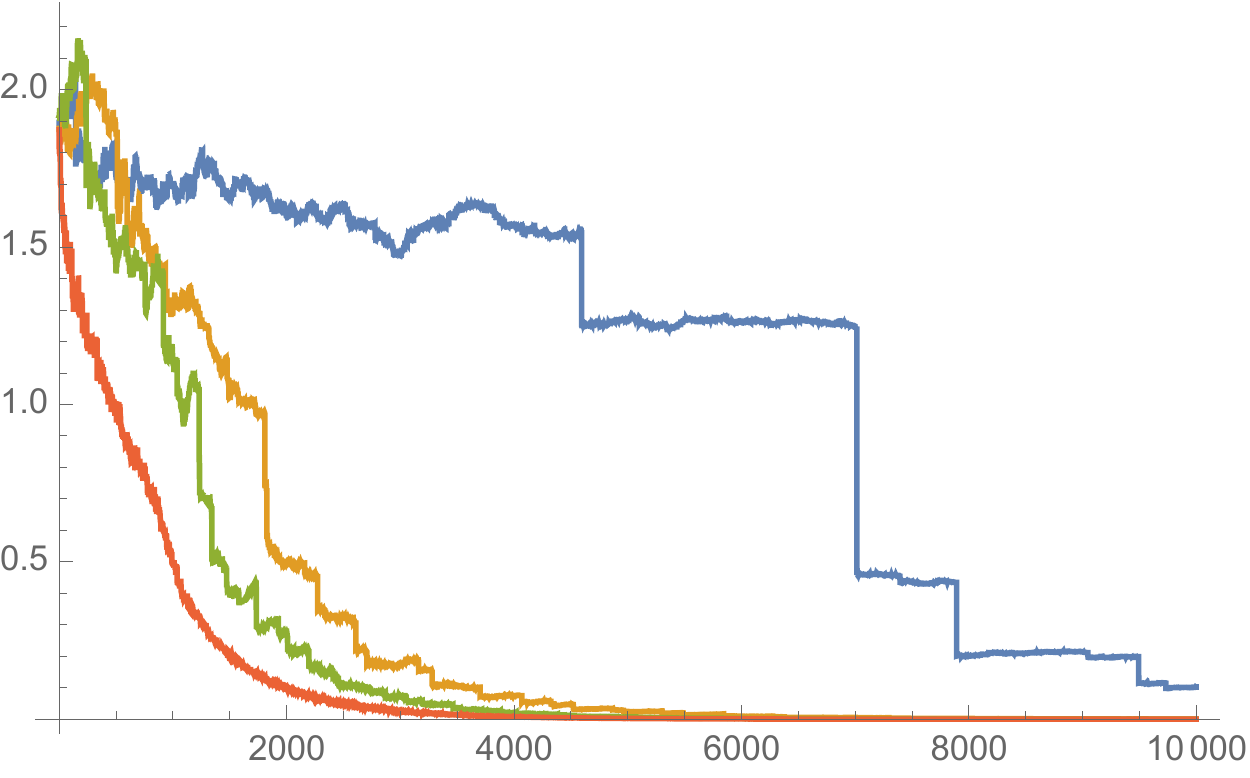}};
\node at (1.2,1.1) {$\|Ax_k - b\|_{\ell^{\infty}}$};
\end{tikzpicture}
\caption{$\|A x_k - b\|_{\ell^{\infty}}$ for the Randomized Kaczmarz method (blue), for $p=1$ (orange), $p=2$ (green) and $p=20$ (red).}
\end{figure}
\end{center}

We observe that the rate of decay is again increasing in $p$.  A priori it is not at all clear that fixing large errors with higher likelihood will necessarily lead to faster decay in the uniform approximation (one may simply create errors somewhere else). This may be a consequence of our using a random matrix. In any case, it would be interesting to have a better understanding of $\|Ax_k - b\|_{\ell^{\infty}}$ both for Randomized Kaczmarz and our variation on it.

\subsection{How good is this bound?} It has been recently shown \cite{stein} that the Randomized Kaczmarz method converges along small singular vectors. More precisely, the Randomized Kaczmarz method satisfies an in interesting identity: if $v_{\ell}$ is a right singular vector of $A$ associated to the singular value $\sigma_{\ell}$, then
$$\mathbb{E} \left\langle x_{k} - x, v_{\ell} \right\rangle = \left(1 - \frac{\sigma_{\ell}^2 }{\|A\|_F^2}\right)^k\left\langle x_0 - x, v_{\ell}\right\rangle.$$
\begin{center}
\begin{figure}[h!]
\begin{tikzpicture}
\node at (0,0) {\includegraphics[width=0.6\textwidth]{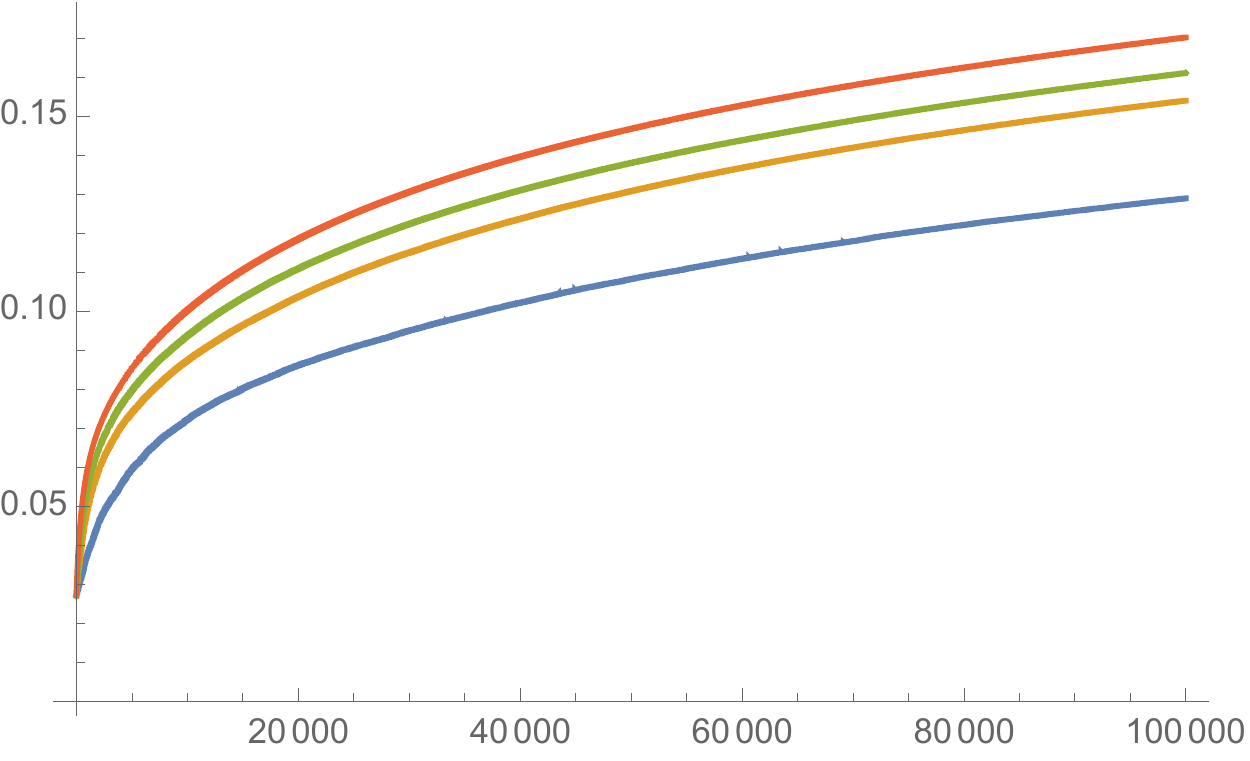}};
\node at (1.5,0) {$\left| \left\langle \frac{x_k - x}{\|x_k - x\|}, v_{1000} \right\rangle \right|$};
\end{tikzpicture}
\caption{The inner product of $x_k - x$ against the smallest singular vector $v_{1000}$ for the Randomized Kaczmarz method (blue), for $p=1$ (orange), $p=2$ (green) and $p=20$ (red).}
\end{figure}
\end{center}
In particular, small singular values induce the slowest decay and, as a consequence, $x_k - x$ tends to have a large component in the space spanned by singular vectors associated to small singular values. It is natural to ask whether the same phenomenon exists for our class of algorithms: it does indeed seem to persist, at least for certain types of matrices, and, in fact, seems to get more pronounced. In Fig. 5, we consider again $A \in \mathbb{R}^{1000 \times 1000}$ by taking each entry to be $a_{ij} \sim \mathcal{N}(0,1)$ and then normalize the rows to $\|a_i\|=1$. We observe that the iterates again seem to concentrate on low-lying singular vectors. 
It would be interesting to understand why that is the case or at least for which types of matrices one could expect this type of phenomenon. We quickly recall the argument from \cite{stein}: assuming w.l.o.g. that we solve $Ax = 0$, that $\|a_i\|=1$ and abbreviating $p_i$ to be the likelihood of choosing the $i-$th equation, we can write
\begin{align*}
\mathbb{E} \left\langle x_{k+1}, v_{\ell} \right\rangle &=  \left\langle x_{k} - \sum_{i=1}^{n}{p_i \left\langle x_k, a_k \right\rangle a_k}, v_{\ell} \right\rangle\\
 &=  \left\langle x_{k}, v_{\ell} \right\rangle - \sum_{i=1}^{n}{p_i\left\langle x_k, a_k \right\rangle \left\langle a_k, v_{\ell} \right\rangle}.
\end{align*}
If the $p_i$ are constant, then the last term is merely 
$$ \sum_{i=1}^{n}{p_i\left\langle x_k, a_k \right\rangle \left\langle a_k, v_{\ell} \right\rangle} = \frac{1}{n} \left\langle Ax_k, A v_{\ell}\right\rangle$$
which then simplifies further since the singular vector $v_{\ell}$ is an eigenvector of $A^TA$. In our case, the $p_i$ are not constant and the argument does not apply. Nonetheless a similar phenomenon seems to be at work -- at least for certain types of matrices. One could also wonder whether the role of small singular vectors is now perhaps played by the minimizer of the functional
$$ J_p(z)  = \frac{\|A z \|^{p+2}_{\ell^{p+2}}}{\|Az \|^p_{\ell^p}\|z\|^2_{2}}$$
which, at least in our Theorem, plays the analogous role that the smallest singular vector plays for the Random Kaczmarz method.

\section{Proofs}
\subsection{Proof of First Part of the Theorem}
\begin{proof} We first prove that
 $$ \mathbb{E} \left\| x_{k+1} - x \right\|_2^2 \leq \left(1 -  \frac{\|A (x-x_k) \|^{p+2}_{\ell^{p+2}}}{\|A(x- x_k) \|^p_{\ell^p}\|x-x_k\|^2_{2}}\right) \|x_k - x\|_2^2$$
which implies the first part of the Theorem. 
 We introduce the error term
$$ e_k = x_k - x.$$
Plugging in, we obtain that if the $i-$th equation is chosen, then
\begin{align*}
x + e_{k+1} &= x_{k+1} = x_k + \frac{ b_i - \left\langle a_i, e_k\right\rangle}{\|a_i\|^2} a_i \\
&=x + e_k + \frac{ b_i - \left\langle a_i, x + e_k\right\rangle}{\|a_i\|^2} a_i \\
&= x + e_k  - \frac{  \left\langle a_i,  e_k\right\rangle}{\|a_i\|^2} a_i + \left( \frac{ b_i - \left\langle a_i, x\right\rangle}{\|a_i\|^2} a_i\right).
\end{align*}
Since $x$ is the exact solution, we have $b_i - \left\langle a_i, x\right\rangle=0$ and
$$ e_{k+1} = e_k - \frac{  \left\langle a_i,  e_k\right\rangle}{\|a_i\|^2} a_i.$$
The Pythagorean Theorem implies that, since $e_{k+1}$ and $a_i$ are orthogonal by design,
$$ \left\|e_{k+1}\right\|^2 = \left\|e_k\right\|^2 - \frac{\left|\left\langle a_i, e_k\right\rangle\right|^2}{\|a_i\|^2}.$$
We note that each row of the matrix is normalized to have row 1 and therefore $\|a_i\|^2 = 1$ and thus
$$ \left\|e_{k+1}\right\|^2 = \left\|e_k\right\|^2 - \left|\left\langle a_i, e_k\right\rangle\right|^2.$$
Let us now assume that
$$ \mathbb{P}(\mbox{we choose equation}~i) =  \frac{\left| \left\langle a_i, x_k \right\rangle - b_i \right|^p}{\|Ax_k - b \|^p_{\ell^p}}.$$
We observe that
$$ \frac{\left| \left\langle a_i, x_k \right\rangle - b_i \right|^p}{\|Ax_k - b \|^p_{\ell^p}} =  \frac{\left| \left\langle a_i, e_k \right\rangle \right|^p}{\|A e_k \|^p_{\ell^p}} $$
and thus, recalling the normalization $\|a_i\|=1$, the expected value of choosing the $i-$th equation with respect to this particular weighting scheme is
\begin{align*}
 \mathbb{E} \left|\left\langle a_i, e_k\right\rangle\right|^2 &= \sum_{i=1}^{m}\left|\left\langle a_i, e_k\right\rangle\right|^2  \frac{\left| \left\langle a_i, e_k \right\rangle \right|^p}{\|A e_k \|^p_{\ell^p}} = \frac{\|A e_k \|^{p+2}_{\ell^{p+2}}}{\|A e_k \|^p_{\ell^p}}.
\end{align*}
This, in turn, implies that
\begin{align*}
 \mathbb{E} \|e_{k+1}\|^2 &= \|e_k\|^2 -  \frac{\|A e_k \|^{p+2}_{\ell^{p+2}}}{\|A e_k \|^p_{\ell^p}} \\
 &= \left( 1 -  \frac{\|A e_k \|^{p+2}_{\ell^{p+2}}}{\|A e_k \|^p_{\ell^p} \|e_k\|_{2}^2} \right) \|e_k\|_{2}^2 \\
  &\leq \left( 1 -  \inf_{x \neq 0} \frac{\|A x \|^{p+2}_{\ell^{p+2}}}{\|A x \|^p_{\ell^p} \|x\|_{2}^2} \right) \|e_k\|_{2}^2.
 \end{align*}
\end{proof}

\subsection{Proof of Second Part of the Theorem}
\begin{proof}
We use H\"older's inequality to argue that for any $0 \neq x \in \mathbb{R}^{m}$
\begin{align*}
\|x\|_{\ell^p}^{p} = \sum_{i=1}^{m}{|x_i|^p} &\leq \left( \sum_{i=1}^{m}{|x_i|^{p+2}}\right)^{\frac{p}{p+2}}  \left( \sum_{i=1}^{m}1\right)^{\frac{2}{p+2}} = \|x\|_{\ell^{p+2}}^{p} m^{\frac{2}{p+2}}.
\end{align*}
This implies, since $Ax \in \mathbb{R}^m$,
$$ \frac{ \|Ax\|_{\ell^{p+2}}^{p+2}}{\|Ax\|_{\ell^p}^p} =\|Ax\|_{\ell^{p+2}}^{2} \frac{ \|Ax\|_{\ell^{p+2}}^{p}}{\|Ax\|_{\ell^p}^p} \geq \frac{\|Ax\|_{\ell^{p+2}}^2}{m^{\frac{2}{p+2}}}.$$
Using H\"older's inequality once more, we have for any $0 \neq x \in \mathbb{R}^{m}$
$$ \|x\|_{\ell^{p+2}} \geq m^{\frac{1}{p+2} - \frac{1}{2}} \|x\|_{2}$$
implying that
$$ \frac{ \|Ax\|_{\ell^{p+2}}^{p+2}}{\|Ax\|_{\ell^p}^p} \geq   \frac{\|Ax\|_{\ell^{p+2}}^2}{m^{\frac{2}{p+2}}} \geq \frac{1}{m} \|Ax\|_{2}^2.$$
Collecting all these estimates, we have that
$$ \frac{\|A x \|^{p+2}_{\ell^{p+2}}}{\|A x \|^p_{\ell^p} \|x\|_{2}^2} \geq \frac{1}{m} \frac{\|Ax\|_{2}^2}{\|x\|_{2}^2}.$$
Recalling the definition of the smallest singular value, or, alternatively, the inverse operator norm, we have
$$ \|Ax\|_{2}^2 \geq \frac{\|x\|_{2}^2}{\| A^{-1}\|_{2}^2}$$
which implies
$$ \frac{\|A x \|^{p+2}_{\ell^{p+2}}}{\|A x \|^p_{\ell^p} \|x\|_{2}^2} \geq \frac{1}{m } \frac{1}{\| A^{-1}\|_{2}^2}.$$
Recall that, using the normalization $\|a_i\|=1$, we have
$$ \|A\|_F^2 = m$$
and we have thus established the desired statement. As for cases of equality, we remark that our use of H\"older's inequality is only sharp if the vector
is constant, which requires $Ax$ to be constant. As for the lower bound on $\|Ax\|_{2}^2,$ it is only sharp if $x$ is the singular
vector associated to the smallest singular value of $A$.
\end{proof}

\end{document}